\documentclass[onecolumn,draftcls,letterpaper,12 pt]{ieeeconf} 

\let\labelindent\relax
\usepackage{times} 
\usepackage{amsmath} 
\usepackage{amssymb,amsfonts}  
\usepackage{latexsym,theorem,epsfig}
\usepackage{graphicx} 
\usepackage{dsfont}
\usepackage{mathtools}
\usepackage{subfigure}
\usepackage{lscape}
\usepackage{float}
\usepackage{enumerate}
 \usepackage[usenames,dvipsnames]{pstricks}
 \usepackage{epsfig}
 \usepackage{pst-grad} 
 \usepackage{pst-plot} 
 \usepackage[percent]{overpic}
 \usepackage{enumitem}  

\usepackage{cite, url}
\usepackage{color,hyperref}
\definecolor{darkblue}{rgb}{0,0,1}
\hypersetup{colorlinks,breaklinks,
	linkcolor=darkblue,urlcolor=darkblue,anchorcolor=darkblue,citecolor=darkblue}
 
\newtheorem{theorem}{Theorem}

\newtheorem{assumption}{Assumption}

{\theorembodyfont{\rmfamily} }
\newtheorem{lemma}[theorem]{Lemma}

\newcommand{\vx}[0]{\ensuremath{\boldsymbol{x}}}
\newcommand{\vy}[0]{\ensuremath{\boldsymbol{y}}}
\newcommand{\vw}[0]{\ensuremath{\boldsymbol{w}}}
\newcommand{\bs}[0]{\ensuremath{\boldsymbol{s}}}
\newcommand{\bq}[0]{\ensuremath{\boldsymbol{q}}}
\newcommand{\bz}[0]{\ensuremath{\boldsymbol{z}}}
\newcommand{\Lap}{\textbf{{\L}}}
\newcommand{\one}{\mathbf{1}}
\allowdisplaybreaks
\IEEEoverridecommandlockouts
\title{Geometrically Convergent Distributed Optimization \\with Uncoordinated Step-Sizes}

\author{Angelia Nedi\'{c}, Alex Olshevsky, Wei Shi, and C\'{e}sar A. Uribe
\thanks{A. Nedi\'{c} (\textit{angelianedich@gmail.com}) is with the ECEE Department, Arizona State University. A. Olshevsky and W. Shi (\{alexols,wilburs\}@bu.edu) are with the ECE Deparment, Boston University. C.A. Uribe (\textit{cauribe2@illinois.edu}) is with the Coordinated Science Laboratory, University of Illinois. 
	This research is supported partially by the National Science Foundation under
	grants CNS 15-44953 and AFOSR FA-95501510394}  
}

\begin{document}
\maketitle
\begin{abstract}
A recent algorithmic family for distributed optimization, DIGing's, have been shown to have geometric convergence over time-varying undirected/directed graphs \cite{Nedic2016}. Nevertheless, an identical step-size for all agents is needed. In this paper, we study the convergence rates of the Adapt-Then-Combine (ATC) variation of the DIGing algorithm under uncoordinated step-sizes. We show that the ATC variation of DIGing algorithm converges geometrically fast even if the step-sizes are different among the agents. In addition, our analysis implies that the ATC structure can accelerate convergence compared to the distributed gradient descent (DGD) structure which has been used in the original DIGing algorithm.   
\end{abstract}

\section{Introduction}

Recent advances in networked and distributed systems require the development of scalable algorithms that take into account the decentralized nature of the problem and communication constraints. Formation control \cite{Ren2006,Ren2007}, distributed spectrum sensing \cite{Bazerque2010}, statistical inference and learning \cite{lee2012,Rabbat2004,Nedic2015_2,ned14} are among some areas of application of such algorithms.

The problem of optimal performance of a number of such distributed systems can be modeled as optimization problems where the objective function is the aggregation of local private information distributed throughout the system. 

This paper focuses on problems of the form
\begin{equation}
	\begin{array}{c}\label{eq:F}
	\min\limits_{x\in\mathbb{R}^p}~f(x)=\frac{1}{n}\sum\limits_{i=1}^n
	f^i(x),
\end{array}
\end{equation}
where each function $f^i: \mathbb{R}^p\rightarrow \mathbb{R}$ is held privately by agent $i$ to encode the agent's objective function, e.g. private data. Moreover, the complete systems seeks to solve the joint problem by exchanging information over a network. Such network might correspond to privacy settings or communication constraints.  

Several algorithms have been proposed for the solution of problems of the form \eqref{eq:F} since the 1980s \cite{tsi84,ber89}. Initial approaches for general and possibly time-varying graphs were based in distributed sub-gradients with extensions to handle stochasticity and asynchronous updates \cite{ned09b,Nedic2011,Ram2010}. Such algorithms are flexible for the class of functions and graphs they can handle but are considerably slow. Even for strongly convex functions a diminishing step-size is required which hinders the possibility of linear rates \cite{Duchi2012,Zhu2012,Yuan2012}. Recent studies have achieved linear convergence rates for strongly convex function \cite{Shi2014,Shi2015,Mokhtari2015,Nedic2016,Xi2015,Zeng2015,Qu2016}. Nonetheless, these methods require a careful selection of the step-sizes. 

Recently in \cite{Xu2015,Xu2016}, the authors utilize the Adapt-Then-Combine strategy\footnote{The readers are referred to reference \cite{Sayed2013} for more discussion on different strategies of information diffusion in networks.} to develop an augmented version of the distributed gradient method for distributed optimization over time-invariant graphs. This algorithm is shown to converge for convex smooth objective functions for sufficiently small constant step-size. Moreover, no coordination on the step-sizes are needed. Additionally similar structures of the dynamic average consensus have been explored for more general classes of non-convex functions \cite{Zhu2010}. For non-convex problems the work in \cite{Lorenzo2016,Lorenzo2015,Lorenzo2016b} develops a large class of distributed algorithms by utilizing varios ``function-surrogate modules" thus providing a great flexibility in its use and rendering a new class of algorithms that subsumes many of the existing distributed algorithms. The authors in \cite{Lorenzo2015,Xu2015} simultaneously proposed methods that track the gradient averages.

In this paper we study the Adapt-Then-Combine Distributed Inexact Gradient Tracking (ATC-DIGing) algorithm for the solution of the optimization problem \eqref{eq:F}. Specifically, we show that geometric convergence rates\footnote{Suppose that a sequence $\{x_k\}$ converges to $x^*$ in some norm $\|\cdot\|$. We say that the convergence is R-linear (Geometric) if there exist $\lambda\in(0,1)$ and some positive constant $C$ such that $\|x_k-x^*\|\leq C\lambda^k$ for all $k$. This rate is often referred to as global to be distinguished from the case when the given relations are valid for some sufficiently large indices $k$.} can still be obtained in the studied algorithm for uncoordinated step-sizes. Moreover, under specific conditions the ATC-DIGing algorithm can use step-sizes as large as the centralized case, which improves the stability region of the ATC-structure over the Distributed Gradient Descent (DGD) structure used in the original DIGing algorithm.

This paper is organized as follows. Section \ref{result} presents some preliminary definitions, the proposed algorithm and the main result of this paper. Section \ref{proof} shows the analysis and proof of the proposed algorithm. Section \ref{discussion} discusses the implications of the results of some general remarks and comments about the contributions. Finally Section \ref{conclusions} presents some conclusions and future work.

\textbf{\textit{Notation.}} Each agent $i$ holds a \emph{local copy} of the variable $x$ of the problem in \eqref{eq:F}, which is denoted by $x^i\in\mathbb{R}^p$; its value at iteration/time $k$ is denoted by $x^i_k$. In general per agent information will be represented by superscripts with letter $i$ or $j$ and time indices by subscripts with the letter $k$. We stack the raw version of all $x^i$ into a single matrix $\vx$ such that ${{\vx \in \mathbb{R}^{n\times p}}}$, while its corresponding $i$-th row is denoted by $\left[ \vx\right]_i = (x^i)'$.  We introduce an aggregate objective function of the local variables: $\left[ \boldsymbol{f}(\boldsymbol{x})\right] _i \triangleq f^i(x^i)$, where its gradient is a matrix whose $i$-th row is defined as $\left[ \nabla \boldsymbol{f} (\vx) \right] _i = \nabla f^i(x^i)'$. We say that $\vx$ is \emph{consensual} if all of its rows are identical, i.e., $x^1=x^2=\cdots= x^n$. Furthermore, we let $\mathbf{1}$ denote a column vector with all entries equal to one (its size is to be understood from the context). The bar denotes averages, e.g. $\bar{\vx}\triangleq1 / n \mathbf{1}' \vx$, and its consensus violation is denotes as $\check{\vx}\triangleq\Lap\vx$, where $\Lap=I- 1/ n\one\one'$. We use $\|\vx\|_{\Lap}$ to denote the $\Lap$ weighted (semi)-norm, that is, $\|\vx\|_{\Lap}=\sqrt{\langle \vx,\Lap \vx \rangle}$. Note that since $\Lap=\Lap'\Lap$, we always have $\|\vx\|_\Lap=\|\Lap \vx\|_F$, where $\| \cdot\|_F$ stands for the Frobenius norm. For a tuple of matrices $\vx,\vy,\vw$ with $\vx=\vy+\vw$, in view of the definition of the consensus violation, it holds that
\begin{align*}
\|\check \vx\|_F = \|\vx\|_{\Lap} = \|\vy+\vw\|_{\Lap} \leq \|\vy\|_{\Lap}+\|\vw\|_{\Lap} = \|\check \vy\|_F+ \|\check \vw\|_F.
\end{align*}

\section{Definitions, Algorithm and Main Result}\label{result}

The set of agents $V=\{1,2,\ldots,n\}$ interact over a time-invariant connected undirected graph $\mathcal{G}= \{V,E\}$, where $E$ correspond to the edges in the graph. A pair of agents $(j,i)\in E$ indicates that agent $j$ can exchange information with agent $i$. The neighbors of agent $i$ is a set defined as $N^i=\left\{j\big|(j,i)\in E\right\}$. Additionally, there is a nonnegative doubly-stochastic weight matrix $W$, compliant with the graph $\mathcal{G}$, such that if $(j,i)\in E$ then $\left[ W\right]_{ij} > 0 $ otherwise $\left[ W\right]_{ij} = 0$.

Next we are going to formalize the set of assumptions we will use for our results.

\begin{assumption}\label{basic}
	The graph $\mathcal{G}$ is connected and $W$ is doubly stochastic.
\end{assumption}

Assumption \ref{basic} is recurrent in many distributed optimization algorithms. It guarantees some minimum exchange of information between agents and balancedness of such exchanges. Specifically this assumption can be relaxed without much extra work for the case of uniformly connected time-varying directed graphs \cite{Nedic2016}.

\begin{lemma}\label{lemma:mixing_contraction}
	Let Assumption \ref{basic} hold. For any matrix $\vx$ with appropriate dimensions, if $\vx = W\vy$, then we have $\|\vx\|_{\emph{\Lap}} \leq \delta\|\vy\|_{\emph{\Lap}}$ where $\delta$ is a constant less than $1$.
\end{lemma}

Lemma \ref{lemma:mixing_contraction} is standard in the consensus literature. An explicit expression of $\delta$ in terms of $n$ can be found in \cite{Nedic2009_2} if more specific assumptions are made.

We also need the following two assumptions on the objective functions, which are common for deriving linear (geometric) rates of gradient-based algorithms for strongly convex smooth optimization problems.

\begin{assumption}[\textbf{Smoothness}]\label{ass:smooth}Every function $f^i$ is differentiable and has Lipschitz continuous gradients, i.e., there exists a constant $L^i\in(0,+\infty)$ such that
\begin{equation}
\begin{array}{c}\nonumber
\|\nabla f^i(x) - \nabla f^i(y)\|_F \leq L^i \|x- y\|_F\ \text{ for any }x, y \in \mathbb{R}^p.
\end{array}
\end{equation}
\end{assumption}

In Section \ref{proof} we will use $L\triangleq\max_i \{L^i\}$, which is the Lipschitz constant of $\boldsymbol{f}(\vx)$, and ${{\bar{L}\triangleq (1/n)\sum_{i=1}^n L^i}}$ which is the Lipschitz constant of $\nabla f(x)$.

\begin{assumption}[\textbf{Strong convexity}]\label{ass:strongly_convex}
	Every function $f^i$ satisfies
	\begin{equation}
	\begin{array}{c}\nonumber
	f^i(x)\geq f^i(y)+\langle\nabla f^i(y),x-y\rangle+\frac{\mu^i}{2}\|x-y\|_F^2,
	\end{array}
	\end{equation}
	for any $x,y \in \mathbb{R}^p$, where $\mu^i\in[0,+\infty)$. Moreover, at least one $\mu^i$ is nonzero.
\end{assumption}

In the analysis we will use ${{\hat{\mu}\triangleq\max_i\{\mu^i\}}}$ and ${{\bar{\mu}\triangleq1/n\sum_{i=1}^n \mu^i}}$. Assumption \ref{ass:strongly_convex} implies the $\bar{\mu}$-strong convexity of $f(x)$. Under this assumption, the optimal solution to problem \eqref{eq:F} is guaranteed to exist and to be unique since $\bar{\mu}>0$. We note that all the convergence results in our analysis are achieved under Assumption \ref{ass:strongly_convex}. We will also use $\bar{\kappa}\triangleq L/\bar{\mu}$.

With the above definitions and assumptions in place, we now state the ATC-DIGing algorithm in its compact vector form. Each agent will maintain two variables $x^i_{k}$ and $y^i_k$ at each time instant $k$. These variables are updated according to the following rule:
\begin{subequations}\label{uncoord_diging}
\begin{align}
\vx_{k+1} & =W\left( \vx_k-D\vy_k \right) \\
\vy_{k+1}& =W\left(  \vy_k+ \nabla \boldsymbol{f}(\vx_{k+1})- \nabla \boldsymbol{f}(\vx_{k})\right) 
\end{align}
\end{subequations}
where $W$ is a doubly stochastic matrix of weights (to be defined soon) and $D$ is a diagonal matrix where $\left[D \right]_{ii} = \alpha_i$ is the step-size of agent $i$. The initial value $\vx_0$ is arbitrary and $\vy_0 = \nabla \boldsymbol{f}(\vx_0)$. 

Algorithm of the form \eqref{uncoord_diging} have been recently proposed under the name Aug-DGM by \cite{Xu2015,Xu2016}, where the convergence of the algorithm under uncoordinated step-sizes is prove. Our objective will be to study the convergence rate of the algorithm. We will show that such algorithm converges geometrically fast, and we will provide an explicit rate estimate. These contributions are stated in the next theorem, which is the main result of this paper.

\begin{theorem}[\textbf{Explicit geometric rate}]\label{theorem:final_bound}
Let Assumptions \ref{basic}, \ref{ass:smooth} and \ref{ass:strongly_convex} hold. Let the step-size matrix $D$ be such that its largest positive entry $\alpha_{\max} = \max_i \alpha_i$ satisfies the following relation:
\begin{align*}
\alpha_{\max}\in\left(0,\min\left\{\frac{(1-\delta)\left(1-\delta-4\sqrt{3}\bar{\kappa}(1-\kappa_D^{-1})\right)}{10L\delta\sqrt{n}\sqrt{\bar{\kappa}}},\frac{1}{2\bar{L}}\right\}\right),
\end{align*}
where $\kappa_D = \alpha_{\max} / \alpha_{\min}$ is the condition number of the step-size matrix $D$, and $\bar{\kappa} = L / \bar{\mu}$. Then, assuming that the step-size heterogeneity is small enough ($\kappa_D<1+\frac{\lambda-\sigma}{4\sqrt{3}\bar{\kappa}}$), the sequence $\{\vx_k\}$ generated by the ATC-DIGing algorithm with uncoordinated step-sizes converges to the optimal solution $\vx^* = \one (x^*)'$ at a global R-linear (geometric) rate $O(\lambda^k)$ where $\lambda \in (0,1)$ is given by
\begin{align}\label{eq:exp_rate_2}
\lambda=\max\left\{\sqrt{12\bar{\kappa}^2(1-\kappa_D^{-1})^2+10L\delta\sqrt{n}\sqrt{\bar{\kappa}}\alpha_{\max}} +\delta+2\sqrt{3}\bar{\kappa}(1-\kappa_D^{-1}),\sqrt{1-\frac{\alpha_{\max}\bar{\mu}}{3}}\right\}
\end{align}
\end{theorem}

Theorem \ref{theorem:final_bound} provides an explicit convergence rate estimate for the ATC-DIGing algorithm. Such rate might not be tight and better choices in the analysis will shown result in better bounds. 

\section{The Small Gain Theorem for Linear Rates}\label{proof}

To establish the R-linear rate of the algorithm, one of our technical innovations will be to resort to a somewhat unusual version of small gain theorem under a well-chosen metric, whose original version has received an extensive research and been widely applied in control theory \cite{Desoer2009}. We choose to analyze the ATC-DIGing algorithm using the small gain theorem due to its effectiveness in showing geometric rates for other algorithms, e.g. \cite{Nedic2016}. We will give an intuition of the whole analytical approach shortly, after stating the small gain theorem at first.

Let us adopt the notation $\bs^i$ for the infinite sequence $\bs^i=\left(\bs^i_0,\bs^i_1,\bs^i_2,\ldots\right)$ where $\bs^i_k\in\mathbb{R}^{n\times p},\ \forall i$. Furthermore, let us define
\begin{subequations}\label{eq:norm}
\begin{align}
\|\bs^i\|_F^{\lambda, K} & \triangleq \max_{k=0,\ldots,K} \frac{1}{\lambda^k}  \|\bs^i_k\|_F \ \ \\
\|\bs^i\|_F^{\lambda} & \triangleq \sup_{k \geq 0} \frac{1}{\lambda^k} \|\bs^i_k\|_F,
\end{align}
\end{subequations}
where the parameter $\lambda\in(0,1)$ will serve as the linear rate parameter later in our analysis. While $\|\bs^i\|_F^{\lambda, K}$ is always finite, $\|\bs^i\|_F^{\lambda}$ may be infinite. If $n=p=1$, i.e., each $\bs^i_k$ is a scalar, we will just write $|\bs^i|^{\lambda, K}$ and $|\bs^i|^{\lambda}$ for these quantities. Intuitively, $\|\bs^i\|_F^{\lambda, K}$ is a weighted ``ergodic norm'' of $\bs^i$. Noticing that the weight $1 / \lambda^k$ is exponentially growing with respect to $k$, if we can show that $\|\bs^i\|_F^{\lambda}$ is bounded, then it
would imply that $\|\bs^i_k\|_F\to0$ geometrically fast. This ergodic definition enables us to give analysis to those algorithms which do not converge Q-linearly. Next we will state the small gain theorem which gives a sufficient condition to for the boundedness of $\|\bs^i\|_F^{\lambda}$. The theorem is a basic result in control systems and a detailed discussion about its result can be found in \cite{Desoer2009}.

\begin{theorem}[\textbf{The small gain theorem}]\label{theorem:small} Suppose that $\{\bs^1, \ldots, \bs^m\}$ is a set of sequences such that for all positive integers $K$ and for each $i=1, \ldots, m$then
	\begin{align}\label{eq:smallgainrecur}
	\|\bs^{(i~{\rm mod}~m)+1} \|_F^{\lambda, K} \leq \gamma_i  \|\bs^{i} \|_F^{\lambda, K} + \omega_i,
	\end{align}
	where the constants (gains) $\gamma_1, \ldots, \gamma_m$ are nonnegative and satisfy $\gamma_1 \gamma_2 \cdots \gamma_m < 1$. Then
	\begin{align}
	\|\bs^1\|_F^{\lambda}&\leq \left( 1 / 1 - \prod_{i=1}^{m}\gamma_i \right) \sum_{i=1}^{m}\omega_i\prod_{j=i+1}^{m}\gamma_j.
	\end{align}
\end{theorem}
For simplicity of exposition we will denote the bound relation in \eqref{eq:smallgainrecur} as an arrow $\bs^i\rightarrow\bs^{(i \mod m)+1}$. Clearly, the small gain theorem involves a cycle $\bs^1\rightarrow\bs^2\rightarrow\cdots\rightarrow\bs^m\rightarrow\bs^1$. Due to this cyclic structure  similar bounds hold for $\|\bs^i\|_F^{\lambda},\ \forall i$.

\begin{lemma}[\textbf{Bounded norm $\Rightarrow$ R-linear rate}]\label{lemma:bound_R-linear}
	For any matrix sequence $\bs^i$, if $\|\bs^i\|_F^{\lambda}$ is bounded, then $\|\bs^i_k\|_F$ converges at a global R-linear (geometric) rate $O(\lambda^k)$.
\end{lemma}

Before summarizing our main proof idea, let us define some quantities which we will use frequently in our analysis. We define $\vx^*\triangleq\one(x^*)'$ where $x^*$ is the optimal solution of \eqref{eq:F}. Also, define
\begin{align}
\bq_k\triangleq\vx_k - \vx^*\ \text{ for any }\ k=0,1,\ldots,
\end{align}
which  is the optimality residual of the iterates $\vx_k$ (at the $k$-th iteration). Moreover, let us adopt the notation
\begin{align}
\bz_k\triangleq \nabla \boldsymbol{f}(\vx_k) - \nabla \boldsymbol{f}(\vx_{k-1})\ \text{ for any }\ k=1,2,\ldots,
\end{align}
and with the convention that $\bz_0\triangleq 0$.

\begin{figure}[H]
	\centering
	\begin{overpic}[width=0.4\textwidth]{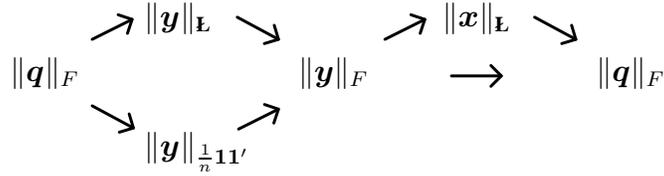}
		\put(-15,12){$\|\bq\|_F$}
		\put(43,12){$\|\vy\|_F$}
		\put(12,23){$\|\vy\|_\Lap$}
		\put(72,23){$\|\vx\|_\Lap$}
		\put(12,-3){$\|\vy\|_{\frac{1}{n}\one\one'}$}
		\put(103,12){$\|\bq\|_F$}
	\end{overpic} 
	\caption{Bound relations between variables in the ATC-DIGing algorithm}
	\label{eq:cycle_alg1}
\end{figure}

We will apply the small gain theorem with the $\|\cdot\|_F^{\lambda, K}$ metric and a right choice of $\lambda < 1$ around the cycle of arrows shown in Figure \ref{eq:cycle_alg1}.

After the establishment of each arrow/relation, we will apply the small gain theorem. Specifically we will use the sequences $\{\|\bq\|_F,\|\vy\|_F,,\|\vy\|_\Lap,\|\vx\|_\Lap,\|\vy\|_{\frac{1}{n}\one\one'},\|\bq\|_F\}$ to show they are bounded and hence conclude that all quantities in the ``circle of arrows'' decay at an R-linear rate $O(\lambda^k)$. 

Note that to apply the small gain theorem, we would need to have gains ($\gamma_i$) that multiply to less than one. This is achieved by choosing an appropriate step-size matrix $D$.

The next lemma presents the establishment of each arrow/relation in the sketch in Fig. \ref{eq:cycle_alg1}.

\begin{lemma}\label{main_lemma}
	Let Assumptions \ref{basic}, \ref{ass:smooth} and \ref{ass:strongly_convex} hold and let $\delta$ be as given in Lemma \ref{lemma:mixing_contraction}. Also, let $\lambda$ be such that $\delta < \lambda < 1$. Then, we have for all $K=0,1,\ldots$,
	\begin{enumerate}[label=(\roman*)]
		\item $ \|\bq\|_F\to\|\vy\|_\Lap:\|\vy\|_{\Lap}^{\lambda, K}  \leq \gamma_{11} \|\bq\|_F^{\lambda, K}+\omega_{11}$
		\item $\|\bq\|_F\to\|\vy\|_{\frac{1}{n}\one\one'}:\|\vy\|_{\frac{1}{n}\one\one'}^{\lambda, K} \leq \gamma_{12} \|\bq\|_F^{\lambda, K} + \omega_{12}$
		\item ${{\{\|\vy\|_{\Lap},\|\vy\|_{\frac{1}{n}\one\one'}\}\to\|\vy\|_F :\|\vy\|_F^{\lambda,K}=\|\vy\|_{\Lap}^{\lambda,K}+\|\vy\|_{\frac{1}{n}\one\one'}^{\lambda,K}}}$
		\item $ \|\vy\|_F\to\|\vx\|_\Lap:\|\vx\|_{\Lap}^{\lambda,K}\leq\gamma_2\|\vy\|_F^{\lambda,K}+\omega_2$
	\end{enumerate}
	where
	\begin{align*}
	\gamma_{11}=\frac{(\lambda+1)\delta L}{\lambda-\delta} \qquad\gamma_{12}=L \qquad \gamma_2=\frac{\delta\alpha_{\max}}{\lambda-\delta}
	\end{align*}
	and $\omega_{11},\omega_{12},\omega_2<+\infty$
\end{lemma} 

Lemma \ref{main_lemma} provides a subset of the required relations necessary for the application of the small gain theorem. Relations $\{\|\vx\|_\Lap, \|\vy\|_F\} \to \|\bq\|_F$ remains to be addressed. For this, we need an interlude on gradient descent with errors in the gradient. Since this part is relatively independent from the preceding development, we provide it in the next subsection.

\subsection{The Inexact Gradient Descent on a Sum of Strongly Convex Functions}

In this subsection, we consider the basic (centralized) first-order method for problem \eqref{eq:F} under inexact first-order oracle.  To distinguish from the notation used for our distributed optimization problem/algorithm/analysis, let us make some definitions that are only used in this subsection. Problem \eqref{eq:F} is restated as follows with different notation,
\[\label{eq:F_diff}
\min\limits_{p\in\mathbb{R}^d} g(p)=\frac{1}{n}\sum\limits_{i=1}^n g^i(p),
\]
where all $g_i$'s satisfy Assumptions \ref{ass:smooth} and \ref{ass:strongly_convex} with $f^i$ being replaced by $g_i$. Let us consider the inexact gradient descent (IGD) on the function $g$:
\begin{align}\label{eq:IGM}
p_{k+1} = p_k - \theta \frac{1}{n} \sum_{i=1}^n \nabla g^i(s^i_k)+e_k,
\end{align}
where $\theta$ is the step-size and $e_k$ is an additive noise. Let $p^*$ be the global minimum of $g$, and define
\[
r_k\triangleq \|p_k - p^*\|_F\text{ for any }k=0,1,\ldots.
\]

The main lemma of this subsection is stated next; it is basically obtained by following the ideas in \cite{Devolder2013,Nedic2016}.

\begin{lemma}[\textbf{The error bound on the IGD}]\label{lemma:graderror}
	Suppose that
	\begin{align}\label{eq:IGM_lambda}
	\sqrt{1-\frac{\theta\bar{\mu}\beta}{2(\beta+1)}}\leq\lambda<1\quad\text{and}\quad \theta\leq\frac{1}{(1+\eta)\bar{L}},
	\end{align}
	where $\beta\geq2$ and $\eta>0$. Let Assumptions \ref{ass:smooth} and \ref{ass:strongly_convex} hold for all $g_i$'s. Then, for a set of subsequences $\{s^1_k,s^2_k,\ldots,s^n_k\}$, the tuple sequence $\{r_k, p_k \}$ generated by the inexact gradient method \eqref{eq:IGM} obeys
	\begin{align}\label{eq:IGM_ineq}
	|r|^{\lambda,K} &\leq  (\lambda\sqrt{n})^{-1}\left(\sqrt{\frac{L(1+\eta)}{\bar{\mu}\eta}+\frac{\hat{\mu}}{\bar{\mu}}\beta}\right)\sum\limits_{i=1}^n\|p-s^i\|_F^{\lambda,K} +2r_0+\frac{\sqrt{3-\theta\bar{\mu}}}{\lambda\theta\bar{\mu}}\|e\|_F^{\lambda,K}
	\end{align}
\end{lemma}

Now we prove the last arrow of our proof sketch [cf. \eqref{eq:cycle_alg1}] in the following lemma. Its establishment will use the error bound on the IGD of Lemma \ref{lemma:graderror}, as a key ingredient.

\begin{lemma}[ $\{\|\vx\|_\Lap, \|\vy\|_F\} \to \|\bq\|_F$]\label{lemma:last_arrow}
	Let Assumptions \ref{basic}, \ref{ass:smooth}, and \ref{ass:strongly_convex} hold. In addition, suppose that the parameters $\alpha$ and $\lambda$ are such that
	\[\label{eq:last_arrow_lambda}
	\sqrt{1-\frac{\alpha\bar{\mu}\beta}{2(\beta+1)}}\leq\lambda<1\quad\text{and}\quad\alpha\leq\frac{1}{(1+\eta)\bar{L}},
	\]
	where $\beta\geq2$ and $\eta>0$ are some tunable parameters. Then, we have
	\begin{align}\label{eq:last_arrow_proof0}
	\|\bq\|_F^{\lambda, K}
	& \leq 
		\left(1+\frac{\sqrt{n}}{\lambda}\left(\sqrt{\frac{L(1+\eta)}{\bar{\mu}\eta}+\frac{\hat{\mu}}{\bar{\mu}}\beta}\right)\right)\|\vx\|_{\Lap}^{\lambda, K}+B_{\vy}\|\vy\|_F^{\lambda,K} +2\sqrt{n}\|\bar{x}_0' -  x^*\|_F,
	\end{align}
	where the constant $B_{\vy}=\frac{\sqrt{3-\alpha_{\max}\bar{\mu}}}{\lambda\bar{\mu}}\left(1-\kappa_D^{-1}\right)$ if $\alpha=\alpha_{\max}$; $B_{\vy}=\frac{\sqrt{3-\bar{\alpha}\bar{\mu}}}{\sqrt{n}\lambda\bar{\mu}\bar{\alpha}}\sqrt{\sum_{i=1}^n\left(\alpha_i-\bar{\alpha}\right)^2}$ if $\alpha=\bar{\alpha}$.
\end{lemma}

\subsection{Proof of Main Result}

We are now ready to show the proof of our main result in Theorem \ref{theorem:final_bound}.

\begin{proof}[Theorem \ref{theorem:final_bound}]
	We will use the small gain Theorem \ref{theorem:small}, together with Lemma \ref{main_lemma} and Lemma \ref{lemma:last_arrow}, to show that $\|\bq\|_F^{\lambda}$ is bounded. Therefore, we need $(\gamma_{11}+\gamma_{12})(\gamma_2\gamma_{31}+\gamma_{32})<1$, that is,
	{
	\begin{align}\label{eq:final_1}
	\left(\frac{\delta\alpha_{\max}}{\lambda-\delta}\left(1+\frac{\sqrt{n}}{\lambda}\sqrt{\frac{L(1+\eta)}{\bar{\mu}\eta}+\frac{\hat{\mu}}{\bar{\mu}}\beta}\right)  +\frac{\sqrt{3-\alpha_{\max}\bar{\mu}}}{\lambda\bar{\mu}}(1-\kappa_D^{-1})\right) 
	 \left(\frac{(\lambda+1)\delta L}{\lambda-\delta}+L\right) <1,
	\end{align}
	}
	where $\beta\geq2$ and $\eta>0$, along with other restrictions on parameters that appear in Lemmas \ref{main_lemma} and \ref{lemma:last_arrow}.

	To obtain some concise though probably loose bound on the convergence rate, we next use Lemma \ref{main_lemma} with some specific values for the parameters $\beta$ and $\eta$, which yields the desired result. Specifically, let $\beta=2L/\hat{\mu}$ and $\eta=1$. By further using $0.5\leq\lambda<1$, it yields from \eqref{eq:final_1} that
	\begin{align}\label{eq:final_proof1}
	\alpha_{\max}\leq\frac{(\lambda-\delta)\left(\lambda-\delta-4\sqrt{3}\bar{\kappa}(1-\kappa_D^{-1})\right)}{10L\delta\sqrt{n}\sqrt{\bar{\kappa}}},
	\end{align}
	where we require/assume $\kappa_D<1+\frac{\lambda-\delta}{4\sqrt{3}\bar{\kappa}}$ so to have $\delta+4\sqrt{3}\bar{\kappa}(1-\kappa_D^{-1})<\lambda$ to ensure the non-emptiness of \eqref{eq:final_1} (this way the right-hand-side of \eqref{eq:final_proof1} is always positive). \eqref{eq:final_proof1} further implies that 
	{
	\begin{align}\label{eq:final_proof1_2}
	\delta+2\sqrt{3}\bar{\kappa}(1-\kappa_D^{-1})+\sqrt{12\bar{\kappa}^2(1-\kappa_D^{-1})^2+10L\delta\sqrt{n}\sqrt{\bar{\kappa}}\alpha_{\max}}\leq\lambda.
	\end{align}
	}
	Solving for $\alpha_{\max}$ when $\lambda < 1$ gives
	\begin{equation}\label{eq:final_proof1_3}
	\alpha_{\max}<\frac{(1-\delta)\left(1-\delta-4\sqrt{3}\bar{\kappa}(1-\kappa_D^{-1})\right)}{10L\delta\sqrt{n}\sqrt{\bar{\kappa}}}.
	\end{equation}
	Meanwhile, considering that $\beta/(\beta+1)\geq2/3$ we have that
	\begin{equation}\label{eq:final_proof2}
	\sqrt{1-\frac{\alpha_{\max}\bar{\mu}}{3}}\leq\lambda<1.
	\end{equation}
	Aggregating the multiple conditions for $\alpha_{\max}$ and $\lambda$ provides the desired result.
\end{proof}

\section{Discussion}\label{discussion}

Other possible choices of $\beta$, $\eta$, $\alpha_{\max}$, and $\lambda$ exist and may give tighter bounds but here we only aim to provide an explicit estimate of the convergence rate.

If $\kappa_D=1$, i.e., the multiple agents use an identical step-size, we can use a small-gain theorem sketch that is similar to the one in reference \cite{Nedic2016} to obtain the geometric rate which is tighter than the result in this paper. Specifically, in this case, to reach $\varepsilon$-accuracy, the number of iterations needed by DIGing is at the order of $O\left(\frac{n^{0.5}\bar{\kappa}^{1.5}}{(1-\delta)^2}\ln\frac{1}{\varepsilon}\right)$ while that by ATC-DIGing is $O\left(\max\left\{\frac{\delta^2n^{0.5}\bar{\kappa}^{1.5}+1}{(1-\delta)^2},\bar{\kappa}\right\}\ln\frac{1}{\varepsilon}\right)$. See appendix \ref{app:comparison} for more detailed explanation. This comparison of the algorithms shows that ATC-DIGing has faster convergence rate and it is less sensitive to the condition number $\bar{\kappa}$, especially when $\delta$ is small (the network is well-connected). This implies that in the DIGing-family, we should use the ATC structure as much as possible.

Here one of our goals is to demonstrate that the geometric rate can still be obtained even if we use uncoordinated step-sizes. Compared to the case of coordinated (identical) step-size in \cite{Nedic2016}, to allow uncoordinated step-size we have to demonstrate that $\|\vy^k\|_F$ decays geometrically fast instead of only $\|\vy^k\|_\Lap$ does so. Thus more steps in the small-gain theorem sketch [c.f. \eqref{eq:cycle_alg1}] is needed and worse bound on the rate is derived. 

Considering the bounds on Theorem \ref{theorem:final_bound}, there is a trade-off between the tolerance of step-size heterogeneity ($\kappa_D$) and the achievable largest step-size ($\alpha_{\max}$). In addition, Theorem \ref{theorem:final_bound} says that when the graph is well-connected ($\delta$ is small) enough and heterogeneity ($\kappa_D$) is small enough, a step-size as large as $\frac{1}{2\bar{L}}$ can be utilized and the corresponding convergence rate can be as fast as $\lambda=\sqrt{1-\bar{\mu}/(6\bar{L})}$.

To make the paper concise, we analyse ATC-DIGing under a rather simple network setting, i.e., time-invariant undirected graph. But it can be expected that the idea of analysis in this paper can be extended to dealing with Push-DIGing and the other possible variants of DIGing even under the setting of time-varying directed graphs. Thus in the following numerical test, we shall conduct the experiments under tougher situation for the DIGing families.

\section{Numerical test}
In this section, we use numerical experiments to demonstrate the performance of DIGing family under uncoordinated step-sizes. The problem we are solving is decentralized Huber minimization over time-varying undirected graphs. The experiments settings including data/graph generation are the same as those in section 6 of reference \cite{Nedic2016} except that at each iteration and agent, we perturb the base step-size $\alpha^\ominus$ by a random variable $\zeta$ satisfying the uniform distribution over interval $(0.5,1.5)$. In other words, at iteration $k$, agent $i$ uses step-size $\alpha_k^i=\alpha^\ominus\zeta_k^i$ where $\zeta_k^i$ is the random variable generated over agent $i$ at time $k$. Monte Carlo simulation shows that such step-size sequences in the current experiment ($n=12$) results a heterogeneity at mean $\mathbb{E}\{\kappa_D\}\approx2.5$. 

\begin{figure}[ht]
	\begin{center}
		\includegraphics[height=8cm]{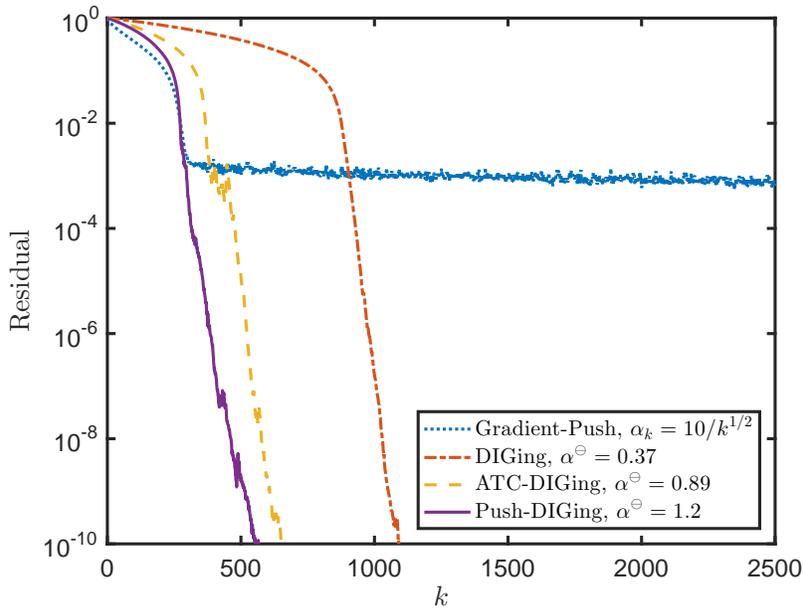}
		\caption{Plot of residuals $\frac{\|\vx^k-\vx^*\|_\mathrm{F}}{\|\vx^0-\vx^*\|_\mathrm{F}}$ for a time-varying undirected graph sequence. The base step-sizes $\alpha^\ominus$'s are set to be the same as the constant step-sizes used in the left sub-figure of Fig. 2 from reference \cite{Nedic2016}.
		}\label{eps:RLS_TV}
	\end{center}
\end{figure}

Numerical results are illustrated in Fig. \ref{eps:RLS_TV}. It shows that under uncoordinated step-sizes, DIGing families still converge geometrically fast.

\section{Conclusions and Future Work}\label{conclusions}

We have shown that the ATC-DIGing algorithm for the distributed optimization problem \eqref{eq:F} converges geometrically to the optimal solution even if all agents have constant uncoordinated step-sizes.  We also provide explicit estimation for its convergence rates. Convergence and rates are derived using the small gain theorem. Nevertheless no claims about tightness of this estimates are given. Under specific conditions the obtained rate shows that ATC-DIGing is less sensitive to the problem parameters that the DGD. Future work should consider extensions to time-varying directed graphs and explore tightness of the rates.  

\bibliographystyle{IEEEtran} 

\bibliography{IEEEfull,bayes_cons_2,document}
\appendices

\section{Proof of Lemma \ref{main_lemma}}

\subsection{(i) $ \|\bq\|_F\to\|\vy\|_\Lap$}
\begin{proof}	
		By the Lipschitz continuity of $\boldsymbol{f}$ (Assumption \ref{ass:smooth}), it is easy to show that for all $K=0,1,\ldots$ and any $\lambda\in (0,1)$,
		\begin{align}\label{eq:first_arrow_proof0}
		\|\bz\|_F^{\lambda, K} \leq L \left( 1 + \frac{1}{\lambda} \right) \|\bq\|_F^{\lambda, K}.
		\end{align}
		
		From \eqref{uncoord_diging}, using Lemma \ref{lemma:mixing_contraction}, it follows that
		\begin{align*}
		\|\vy_{k+1}\|_\Lap
		& \leq \left\|W\vy_k\right\|_{\Lap} + \|W\bz_{k+1}\|_{\Lap}\\
		& \leq \delta\|\vy_k\|_\Lap + \delta\|\bz_{k+1}\|_F,
		\end{align*}
		and therefore, for all $k=0,1,\ldots$,
		\begin{align*}
		\lambda^{-(k+1)}\|\vy_{k+1}\|_\Lap\leq\frac{\delta}{\lambda}\lambda^{-k} 
		\|\vy_k\|_\Lap + \frac{\delta}{\lambda}\lambda^{-k}\|\bz_{k+1}\|_F.
		\end{align*}
		Taking the maximum over $k=0,\ldots,K-1$ on both sides of the above relation, we obtain
		\begin{align*}
		\|\vy\|_\Lap^{\lambda,K}
		\leq\frac{\delta}{\lambda}\|\vy\|_\Lap^{\lambda,K-1}
		+\delta\|\bz\|_F^{\lambda,K}+\|\vy_0\|_{\Lap}
		\end{align*}
		Hence,
		\begin{align}\label{eq:second_arrow_proof2}
		\|\vy\|_\Lap^{\lambda,K}
		&\leq\frac{(\lambda+1)\delta L}{\lambda-\delta} \|\bq\|_F^{\lambda, K}+\frac{\lambda}{\lambda-\delta}\|\vy(0)\|_{\Lap}.
		\end{align}
		Combining \eqref{eq:first_arrow_proof0} and \eqref{eq:second_arrow_proof2} completes the proof.
\end{proof}

\subsection{(ii) $\|\bq\|_F\to\|\vy\|_{\frac{1}{n}\one\one'}$}	

\begin{proof}
Considering $\one'\vy_k-\one'\nabla \boldsymbol{f}(\vx_k)=\cdots=\one'\vy_0-\one'\nabla \boldsymbol{f}(\vx_0)=0$ and $\one'\nabla \boldsymbol{f}(\vx^*)=0$ we have
\begin{align*}
\|\vy\|_{\frac{1}{n}\one\one'}&=\|\frac{1}{\sqrt{n}}\one'(\nabla \boldsymbol{f}(\vx_k)-\nabla \boldsymbol{f}(\vx^*))\|_F\\
&\leq \|\nabla \boldsymbol{f}(\vx_k)-\nabla \boldsymbol{f}(\vx^*)\|_F\\
&\leq L\|\vx_k-\vx^*\|_F.
\end{align*}
The desired result follows immediately.
\end{proof}
\subsection{(iii) $\{\|\vy\|_{\Lap},\|\vy\|_{\frac{1}{n}\one\one'}\}\to\|\vy\|_F $}	This follows automatically from definition.

\subsection{(iv) $ \|\vy\|_F\to\|\vx\|_\Lap$}

\begin{proof}
From \eqref{uncoord_diging}, using Lemma \ref{lemma:mixing_contraction}, for all $k \geq 0$, it follows that
\begin{align*}
\|\vx_{k+1}\|_\Lap
& \leq \left\|W\vx_k\right\|_{\Lap} + \left\|W D\vy_k\right\|_{\Lap}\\
& \leq \delta\|\vx_k\|_\Lap + \delta \alpha_{\max}\|\vy_k\|_F,
\end{align*}
and therefore, for all $k=0,1,\ldots$,
\begin{equation}\label{eq:second2_arrow_proof1}
\lambda^{-(k+1)}\|\vx_{k+1}\|_\Lap\leq\frac{\delta}{\lambda}\lambda^{-k} 
\|\vx_k\|_\Lap + \frac{\delta\alpha_{\max}}{\lambda}\lambda^{-k}\|\vy_k\|_F
\end{equation}
Taking the maximum over $k=0,\ldots,$ on both sides of \eqref{eq:second2_arrow_proof1}, we obtain
\begin{align*}
\|\vx\|_\Lap^{\lambda,K}
&\leq\frac{\delta}{\lambda}\|\vx\|_\Lap^{\lambda,K-1}+\frac{\delta\alpha_{\max}}{\lambda}\|\vy\|_F^{\lambda,K-1}+\|\vx_0\|_{\Lap}\\
&\leq\frac{\delta}{\lambda}\|\vx\|_\Lap^{\lambda,K}+\frac{\delta\alpha_{\max}}{\lambda}\|\vy\|_F^{\lambda,K}+\|\vx_0\|_{\Lap}
\end{align*}
This completes the proof.
\end{proof}

\section{Proof of Lemma \ref{lemma:graderror}}
\begin{proof}
	By assumptions, for each $i\in\{1,2,\ldots,n\}$ and $k=0,1,\ldots$, we have
	\begin{equation}\label{eq:IGM_proof1}
	g^i(p^*) \geq g^i(s_k^i) + \langle\nabla g^i(s_k^i),p^* - s_k^i\rangle + \frac{\mu^i}{2} \| p^* - s_k^i \|_F^2.
	\end{equation}
	Through using the basic inequality $\|s_k^i - p^*\|_F^2 \geq \frac{\beta}{\beta+1}\|p_k - p^*\|_F^2 - \beta\|p_k - s_k^i\|_F^2$ where $\beta>0$ is a tunable parameter, it follows from \eqref{eq:IGM_proof1} that
	\begin{align*}
	g^i(p^*)
	&\geq g^i(s_k^i) +\langle\nabla g^i(s_k^i),p_k - s_k^i\rangle +  \langle\nabla g^i(s_k^i),p^* - p_k\rangle + \frac{\mu^i}{2} \left( \frac{\beta}{\beta+1}\|p_k - p^*\|_F^2 - \beta\|p_k - s_k^i\|_F^2 \right)
	\end{align*}
	and therefore
	\begin{align}\label{eq:IGM_proof3}
	\langle\nabla g^i(s_k^i), p^* - p_k\rangle
	\leq g^i(p^*) - g^i(s_k^i) - \langle\nabla g^i(s_k^i), p_k - s_k^i\rangle  - \frac{\mu^i\beta}{2(\beta+1)} \|p_k - p^*\|_F^2 + \frac{\mu^i\beta}{2} \|s_k^i - p_k\|_F^2.
	\end{align}
	Averaging \eqref{eq:IGM_proof3} over $i$ through $1$ to $n$ gives
	{
	\begin{align} \label{eq:IGM_lb}
 \frac{1}{n} \sum_{i=1}^n \langle\nabla g^i(s_k^i), p^* - p_k\rangle
	&\leq g(p^*)  - \frac{\bar{\mu}\beta}{2(\beta+1)} \|p_k - p^*\|_F^2 \nonumber \\
	& \ 
	 - \frac{1}{n} \sum_{i=1}^n \left( g^i(s_k^i) + \langle\nabla g^i(s_k^i), p_k - s_k^i\rangle - \frac{\mu^i\beta}{2} \|s_k^i - p_k\|_F^2 \right).
	\end{align}
}

	On the other hand, we also have that for any vector $\Delta$,
	{
	\begin{align*}
	 g^i(p_k + \Delta)
 &=   g^i \left(s_k^i + \Delta + p_k - s_k^i \right) \\
	&  \leq g^i(s_k^i) + \langle\nabla g^i(s_k^i), \Delta  + p_k - s_k^i\rangle  + \frac{L^i}{2} \|\Delta + p_k - s_k^i\|_F^2 \\
	& \leq g^i(s_k^i) + \langle\nabla g^i(s_k^i), \Delta\rangle + \langle\nabla g^i(s_k^i),p_k - s_k^i\rangle\\
	&  \qquad   + \frac{L^i(1+\eta)}{2} \|\Delta\|_F^2  + \frac{L^i(1+\eta)}{2\eta} \|p_k - s_k^i\|_F^2
	\end{align*}
	}
	where $\eta>0$ is some tunable parameter, and therefore
	{
	\begin{align}\label{eq:IGM_proof4}
	- \langle\nabla g^i(s_k^i), \Delta\rangle
	& \leq - g^i(p_k + \Delta) + g^i(s_k^i) + \langle\nabla g^i(s_k^i), p_k - s_k^i\rangle \nonumber\\
	& \qquad    + \frac{L^i(1+\eta)}{2\eta} \|p_k - s_k^i\|_F^2 + \frac{L^i(1+\eta)}{2}\|\Delta\|_F^2.
	\end{align}
	}
	Averaging \eqref{eq:IGM_proof4} over $i$ through $1$ to $n$ gives
	{
	\begin{align}\label{eq:IGM_ub}
	&- \langle\frac{1}{n} \sum_{i=1}^n \nabla g^i(s_k^i), \Delta\rangle
	 \leq  - g(p_k + \Delta)+\frac{\bar{L}(1+\eta)}{2}\|\Delta\|_F^2 \nonumber \\
	& \qquad     + \frac{1}{n} \sum_{i=1}^n \left( g^i(s_k^i) + \langle\nabla f^i(s_k^i), p_k - s_k^i\rangle  + \frac{L^i(1+\eta)}{2\eta} \|p_k - s_k^i\|_F^2 \right).
	\end{align}
	}
	Having \eqref{eq:IGM_lb} and \eqref{eq:IGM_ub} at hand, we are ready to show how $r_{k+1}$ is related to $r_k$. First, plugging $a = p_{k+1} - p^*$ and $b = p_k - p_{k+1}$ into the basic equality $\|a\|_F^2=\|a+b\|_F^2 - 2 \langle a, b\rangle - \|b\|_F^2$ yields
	{
	\begin{align}\label{eq:IGM_proof5}
	r_{k+1}^2
	& =  r_k^2 - 2 \langle p_{k+1} - p^*, p_k - p_{k+1} \rangle - \|p_{k+1} - p_k\|_F^2 \nonumber\\
	& =  r_k^2 - 2 \langle p_{k+1} - p^*, \theta \frac{1}{n} \sum\limits_{i=1}^n \nabla g^i(s^i_k) - e_k \rangle  - \|p_{k+1} - p_k\|_F^2 \nonumber\\
	& =  r_k^2 + 2 \theta \langle \frac{1}{n} \sum\limits_{i=1}^n \nabla g^i(s^i_k), p^* - p^{k}\rangle  - 2 \theta \langle \frac{1}{n} \sum\limits_{i=1}^n \nabla g^i(s^i_k), p_{k+1} - p_k \rangle \nonumber\\
	& \qquad  - \|p_{k+1} - p_k\|_F^2 + 2 \langle p_{k+1} - p^*,  e_k \rangle.
	\end{align}
	}
	Next, in \eqref{eq:IGM_proof5}, we substitute \eqref{eq:IGM_lb} for the second term, and we substitute \eqref{eq:IGM_ub} with $\Delta=p_{k+1} - p_k$ for the third term. Thus, we obtain that
	{
	\begin{align}\label{eq:IGM_proof6}
	 r_{k+1}^2
	& \leq r_k^2 + 2\theta \left(g(p^*) - \frac{1}{n} \sum_{i=1}^n \left( g^i(s_k^i) + \langle\nabla g^i(s_k^i), p_k - s_k^i\rangle- \frac{\mu^i\beta}{2} \|s_k^i - p_k\|_F^2 \right)   \right.  
	\nonumber \\ &
	\left. \qquad \qquad - \frac{\bar{\mu}\beta}{2(\beta+1)} \|p_k - p^*\|_F^2\right) + 2\theta \left( - g(p_{k+1})  \right.  
	\nonumber \\ &
	\qquad \qquad +\left.\frac{1}{n} \sum_{i=1}^n \left(  g^i(s_k^i) + \langle\nabla g^i(s_k^i), p_k - s_k^i\rangle + \frac{L^i(1+\eta)}{2\eta} \|p_k - s_k^i\|_F^2 \right) 
	\right.  
	\nonumber \\ &
	\qquad \qquad \left.+\frac{\bar{L}(1+\eta)}{2}\|p_{k+1}-p_k\|_F^2 \right) - \|p_{k+1} - p_k\|_F^2+ 2 \langle p_{k+1} - p^*,  e_k \rangle\nonumber\\
	& \leq r_k^2 + 2\theta ( g(p^*) - g(p_{k+1}) ) - \frac{\theta\bar{\mu}\beta}{\beta+1} \|p_k- p^*\|_F^2  + \frac{1}{n}\sum_{i=1}^n \left(\frac{\theta L^i(1+\eta)}{\eta}+\theta\mu^i\beta\right)\|p_k
	\nonumber \\ &
	\qquad \qquad - s^i_k\|_F^2 - p^*\|_F^2+\frac{1}{\rho}\|e_k\|_F^2 
- (1-\theta \bar{L}(1+\eta))\|p_{k+1}-p_k\|_F^2 + \rho\|p_{k+1} \nonumber\\
	& \leq \left( 1 - \frac{\theta\bar{\mu}\beta}{\beta+1} \right) r_k^2 - 2\theta ( g(p_{k+1}) - g(p^*) ) + \left(\frac{\theta L(1+\eta)}{\eta}
 +\theta\hat{\mu}\beta\right) \frac{1}{n}\sum_{i=1}^n \|p_k - s^i_k\|_F^2
	\nonumber \\ &
	\qquad \qquad
	  - (1-\theta \bar{L}(1+\eta))\|p_{k+1}-p_k\|_F^2 
	+ \rho_kr_{k+1}^2+\frac{1}{\rho_k}\|e_k\|_F^2,
	\end{align}
	}
	where $\{\rho_k\}$ is a sequence of positive tunable parameters (intuitively since in some scenario the noise term $\|e_k\|_F$ decays to zero eventually, a time-wise/diminishing parameter $\rho_k$ may improve the analytical rate). Define $\epsilon_{k}=\frac{1}{n} \sum_{i=1}^n \|p_k - s_k^i\|_F^2$. By choosing $\theta\leq\frac{1}{(1+\eta)\bar{L}}$ such that $1-\theta \bar{L}(1+\eta)$ in \eqref{eq:IGM_proof6} is nonnegative, we have
	{
	\begin{align}\label{eq:IGM_proof7}
	r_{k+1}^2 &\leq \frac{1}{1-\rho_k}\left( 1 - \frac{\theta\bar{\mu}\beta}{\beta+1} \right) r_k^2 - \frac{2\theta}{1-\rho_k}( g(p_{k+1}) - g(p^*) )+ \frac{1}{1-\rho_k}\left(\frac{\theta L(1+\eta)}{\eta}+\theta\hat{\mu}\beta\right)\epsilon_{k}\nonumber\\
	& \qquad \qquad   +\frac{1}{(1-\rho_k)\rho_k}\|e_k\|_F^2.
	\end{align}
	}
	Let us look into the second and third terms in the right-hand-side of \eqref{eq:IGM_proof7}. Noticing that $\bar{\mu}=(1/n)\sum_{i=1}^n\mu^i$ is a strong convexity constant of $g(p)$, there are two possibilities that could happen at time $k$. Possibility A is that
	\[
	r_{k+1}^2\geq \left(\frac{L(1+\eta)}{\bar{\mu}\eta}+\frac{\hat{\mu}}{\bar{\mu}}\beta\right)\epsilon_k+\frac{1}{\rho_k\theta\bar{\mu}}\|e_k\|_F^2,
	\]
	while possibility B is the opposite, namely that
	\[
	r_{k+1}^2 <  \left(\frac{L(1+\eta)}{\bar{\mu}\eta}+\frac{\hat{\mu}}{\bar{\mu}}\beta\right)\epsilon_k + \frac{1}{\rho_k\theta\bar{\mu}}\|e_k\|_F^2.
	\]
	If possibility A occurs, we have
	{
	\begin{align*}
	2\theta(g(p_{k+1})  - g(p^*)) &  \geq  \theta \bar{\mu} \|p_{k+1}-p^*\|_F^2 \\
	&=\theta\bar{\mu} (r_{k+1})^2  \\
	& \geq \left(\frac{\theta L(1+\eta)}{\eta}+\theta\hat{\mu}\beta\right)\epsilon_k+\frac{1}{\rho_k}\|e_k\|_F^2
	\end{align*}
}
which together with \eqref{eq:IGM_proof7} implies
	\[
	r_{k+1}^2\leq\frac{1}{1-\rho_k}\left( 1 - \frac{\theta\bar{\mu}\beta}{\beta+1} \right)r_k^2.
	\]
	
	Considering both possibilities A and B, it follows that
	\begin{align}\label{eq:IGM_proof8}
	r_{k+1}^2& \leq\max\left\{\frac{1}{1-\rho_k}\left(1-\frac{\theta\bar{\mu}\beta}{\beta+1} \right)r_k^2, \left(\frac{L(1+\eta)}{\bar{\mu}\eta}+\frac{\hat{\mu}}{\bar{\mu}}\beta\right)\epsilon_k+\frac{1}{\rho_k\theta\bar{\mu}}\|e_k\|_F^2\right\}.
	\end{align}
	
	Recursively using the inequality \eqref{eq:IGM_proof8}, we get
	{
	\begin{align}\label{eq:IGM_proof9}
	&\lambda^{-2(k+1)}r_{k+1}^2\leq\max
	\left\{\lambda^{-2(k+1)}\left(\prod\limits_{t=0}^k\frac{1}{1-\rho^t}\right)\times \right.  \nonumber \\
	& \left. \left(1-\frac{\theta\bar{\mu}\beta}{\beta+1} \right)^{k+1}r_0^2,\lambda^{-2(k+1)}\max\limits_{t=0,\ldots,k}\left\{\left(\prod\limits_{s=0}^{t-1}\frac{1}{1-\rho_{k-s}}\right)\times \right.  \right.  \nonumber \\
	& \left. \left.  \left(1-\frac{\theta\bar{\mu}\beta}{\beta+1} \right)^t\left(\left(\frac{L(1+\eta)}{\bar{\mu}\eta}+\frac{\hat{\mu}}{\bar{\mu}}\beta\right)\epsilon_{k-t}+\frac{\|e_{k-t}\|_F^2}{\rho_{k-t}\theta\bar{\mu}}\right)\right\}\right\}.
	\end{align}
}
	Taking square root on both sides of \eqref{eq:IGM_proof9} gives us
	{
	\begin{align}\label{eq:IGM_proof10}
	&\lambda^{-(k+1)}r_{k+1}\leq\lambda^{-(k+1)}\left(\prod\limits_{t=0}^k\sqrt{\frac{1}{1-\rho^t}}\right)\left(1-\frac{\theta\bar{\mu}\beta}{\beta+1}\right)^\frac{k+1}{2}r_0 \nonumber \\
	&
	+\lambda^{-(k+1)}\max\limits_{t=0,\ldots,k}\left\{\left(\prod\limits_{s=0}^{t-1}\sqrt{\frac{1}{1-\rho_{k-s}}}\right)\left(\sqrt{1-\frac{\theta\bar{\mu}\beta}{\beta+1} }\right)^t \times \right.  \nonumber \\
	& \left. \left(\sqrt{\frac{L(1+\eta)}{\bar{\mu}\eta}+\frac{\hat{\mu}}{\bar{\mu}}\beta}\sqrt{\epsilon_{k-t}}+\sqrt{\frac{1}{\rho_{k-t}\theta\bar{\mu}}}\|e_{k-t}\|_F\right)\right\}.
	\end{align}
	}
	It can be seen that
	\vskip -1.5pc
	\begin{align*}
	\sup_{k=0,1,\ldots}(\lambda)^{-(k+1)}\left(\prod\limits_{t=0}^k\sqrt{\frac{1}{1-\rho_t}}\right)\left(\sqrt{1-\frac{\theta\bar{\mu}\beta}{\beta+1}}\right)^{k+1}r_0
	\end{align*} 
	is finite (exists) as long for some large enough $t_0$, $\forall t>t_0$ we have $\frac{1}{\lambda}\sqrt{\frac{1}{(1-\rho_t)}\left(1-\frac{\theta\bar{\mu}\beta}{\beta+1}\right)}$ being no greater than $1$. This can be sufficed by setting $\rho_k$ to some specific sequence that monotonically diminishes to zero and requiring that $\lambda>\sqrt{1-\frac{\theta\bar{\mu}\beta}{\beta+1}}$ simultaneously. This can also be done by simply setting a time-invariant $\rho_k$, that is, $\rho_k=\rho\triangleq\frac{\theta\bar{\mu}\beta}{2(\beta+1)-\theta\bar{\mu}\beta}$, along with $\lambda\geq\sqrt{1-\frac{\theta\bar{\mu}\beta}{2(\beta+1)}}$. However, Noticing that the above two different choices of $\rho_k$ lead to bounds on $\lambda$ that have the same order, in the following for conciseness we use the aforementioned time-invariant choice of $\rho_k=\rho$.	
	
	Choosing $\rho_k=\rho$, $\lambda\geq\sqrt{1-\frac{\theta\bar{\mu}\beta}{2(\beta+1)}}$, and $\beta\geq2$, then, \eqref{eq:IGM_proof10} can be further relaxed to
	\vskip -1.5pc
	{
	\begin{align}\label{eq:IGM_proof10_2}
	&\lambda^{-(k+1)}r_{k+1}\leq r_0
	+\lambda^{-(k+1)}\max\limits_{t=0,\ldots,k}\left\{\left(\sqrt{1-\frac{\theta\bar{\mu}\beta}{2(\beta+1)} }\right)^t \right.  \nonumber \\
	& \left. \left(\sqrt{\frac{L(1+\eta)}{\bar{\mu}\eta}+\frac{\hat{\mu}}{\bar{\mu}}\beta}\sqrt{\epsilon_{k-t}}+\frac{\sqrt{3-\theta\bar{\mu}}}{\theta\bar{\mu}}\|e_{k-t}\|_F\right)\right\}.
	\end{align}
	}
	
	Let us denote $c\triangleq\lambda^{-2}\left(1-\frac{\theta\bar{\mu}\beta}{2(\beta+1)}\right)\leq1$, then from \eqref{eq:IGM_proof10_2} we get
	\begin{align}\label{eq:IGM_proof11}
	&\lambda^{-(k+1)}r^{k+1}
	\leq r_0+\lambda^{-1}\max\limits_{t=0,\ldots,k}\left\{\lambda^{-(k-t)}\sqrt{c}^t \times\right.  \nonumber \\
	& \left. \qquad \left(\sqrt{\frac{L(1+\eta)}{\bar{\mu}\eta}+\frac{\hat{\mu}}{\bar{\mu}}\beta}\sqrt{\epsilon_{k-t}}+\frac{\sqrt{3-\theta\bar{\mu}}}{\theta\bar{\mu}}\|e_{k-t}\|_F\right)\right\} \nonumber\\
	&\leq r_0+\lambda^{-1}\max\limits_{t=0,\ldots,k}\left\{\lambda^{-t}\left(\sqrt{\frac{L(1+\eta)}{\bar{\mu}\eta}+\frac{\hat{\mu}}{\bar{\mu}}\beta}\sqrt{\epsilon_{t}}\right. \right. \nonumber \\
	&
	\left. \left. \qquad +\frac{\sqrt{3-\theta\bar{\mu}}}{\theta\bar{\mu}}\|e_{t}\|_F\right)\right\}.
	\end{align}
	Further observing that
	\[
	\sqrt{\epsilon_k} = \sqrt{\frac{1}{n} \sum_{i=1}^n \|p_k - s_k^i\|_F^2} \leq \frac{1}{\sqrt{n}} \sum\limits_{i=1}^n \|p_k - s_k^i\|_F,
	\]
	and combining it with \eqref{eq:IGM_proof11}, it follows that
	{
	\begin{align}\label{eq:IGM_proof12}
	&\lambda^{-(k+1)}r^{k+1}\leq
 r_0+ +\frac{\sqrt{3-\theta\bar{\mu}}}{\lambda\theta\bar{\mu}}\max\limits_{t=0,\ldots,k}\left\{\lambda^{-t}\|e_{t}\|_F\right\} +\lambda\sqrt{n}^{-1} \times\nonumber \\
	& \qquad\left(\sqrt{\frac{L(1+\eta)}{\bar{\mu}\eta}+\frac{\hat{\mu}}{\bar{\mu}}\beta}\right)\sum\limits_{i=1}^n\max\limits_{t=0,\ldots,k}\left\{\lambda^{-t}\|p_t - s_t^i\|_F\right\}
	\end{align}
	}
	Taking $\max_{k=0,1,\ldots,K-1}$ on both sides of \eqref{eq:IGM_proof12}  results in the \eqref{eq:IGM_ineq}.
\end{proof}

\section{Proof of Lemma \ref{lemma:last_arrow}}

\begin{proof}
		First, let us consider the evolution of ${{\bar{x}_k\triangleq\frac{1}{n}\one'\vx_k}}$. Notice that $\one'\vy_k=\one'\nabla \boldsymbol{f}(\vx_k)$ holds for all $k$, we then have that
		\begin{equation}\label{eq:last_arrow_proof1}
		\bar{x}_{k+1}=\bar{x}_k - \alpha\frac{1}{n}\sum\limits_{i=1}^n\nabla f^i(x^i_k)+\left(\alpha\frac{1}{n}\one'-\frac{1}{n}\one' D\right)\vy_k
		\end{equation}
		where $\alpha$ is some nonnegative constant that we consider as the (centralized) step-size of the gradient descent.
		Applying Lemma \ref{lemma:graderror} to the recursion relation of $\bar{x}$, namely \eqref{eq:last_arrow_proof1}, we obtain
		{
		\begin{align}\label{eq:last_arrow_proof2}
		\|\bar{x}' - x^*\|_F^{\lambda, K}
		&\leq  2\|(\bar{x}_0)' - x^*\|_F +(\lambda\sqrt{n})^{-1}\left(\sqrt{\frac{L(1+\eta)}{\bar{\mu}\eta}+\frac{\hat{\mu}}{\bar{\mu}}\beta}\right)\sum\limits_{i=1}^n \|\bar{x}' - x^i\|_F^{\lambda, K} \nonumber\\
		&    
		\qquad+\frac{\sqrt{3-\alpha\bar{\mu}}}{\lambda\alpha\bar{\mu}}\|\left(\alpha\frac{1}{n}\one'-\frac{1}{n}\one'D\right)\vy\|_F^{\lambda,K} \nonumber \\
		&  \leq 2\|(\bar{x}_0)' - x^*\|_F
		+\lambda^{-1}\left(\sqrt{\frac{L(1+\eta)}{\bar{\mu}\eta}+\frac{\hat{\mu}}{\bar{\mu}}\beta}\right)\|\vx\|_\Lap^{\lambda,K}   +\frac{\sqrt{3-\alpha\bar{\mu}}}{n\lambda\bar{\mu}}\|\one'-\one'\frac{D}{\alpha}\vy\|_F^{\lambda,K}.
		\end{align}
	}
		Since $\bq_k = \Lap\vx_k + \frac{1}{n}\one\one'\vx_k  - \vx^*$, it follows that
		\begin{equation}\label{eq:last_arrow_proof4}
		\|\bq\|_F^{\lambda, K} \leq \|\vx\|_\Lap^{\lambda, K} + \sqrt{n}\|\bar{x} - x^*\|_F^{\lambda, K}.
		\end{equation}
		Substituting \eqref{eq:last_arrow_proof2} into \eqref{eq:last_arrow_proof4} yields
		{
		\begin{align}\label{eq:last_arrow_proof5}
		&\|\bq\|_F^{\lambda, K}
		\leq 2\sqrt{n}\|\bar{x}_0' - x^*\|_F\
		+\frac{\sqrt{3-\alpha\bar{\mu}}}{\sqrt{n}\lambda\bar{\mu}}\|\one'-\one'\frac{D}{\alpha}\vy\|_F^{\lambda,K} +\left(1+\frac{\sqrt{n}}{\lambda}\left(\sqrt{\frac{L(1+\eta)}{\bar{\mu}\eta}+\frac{\hat{\mu}}{\bar{\mu}}\beta}\right)\right)\|\vx\|_{\Lap}^{\lambda, K}
		\end{align}
	}
		If we choose $\alpha=\alpha_{\max}$ with $\kappa_D=\frac{\alpha_{\max}}{\alpha_{\min}}$ being the condition number of $D$; and if we choose $\alpha=\frac{1}{n}\sum^i\alpha^i$, we get the desired bounds.
\end{proof}

\section{Comparison of the complexities of DIGing and ATC-DIGing}\label{app:comparison}
	
	Let us first illustrate the main difference between the rates/complexities derived in the current paper and reference \cite{Nedic2016}. The major difference comes from the different restriction of gains that multiply to less than $1$. Let us assume $\kappa_D=1$. In reference \cite{Nedic2016}, the restriction on DIGing is
	\begin{equation}\label{eq:comparison_1}
	\frac{\lambda+1}{\lambda-\delta}\frac{\alpha L}{\lambda-\delta}(5\sqrt{\bar{\kappa}n})<1.
	\end{equation}
	Following the idea of reference \cite{Nedic2016}, the restriction on ATC-DIGing will be
	\begin{equation}\label{eq:comparison_2}
	\frac{\delta(\lambda+1)}{\lambda-\delta}\frac{\delta\alpha L}{\lambda-\delta}(5\sqrt{\bar{\kappa}n})<1.
	\end{equation}
	In the current paper's scheme, the restriction on ATC-DIGing is
	\begin{equation}\label{eq:comparison_3}
	\left(\frac{\delta(\lambda+1)}{\lambda-\delta}+1\right)\frac{\delta\alpha L}{\lambda-\delta}(5\sqrt{\bar{\kappa}n})<1.
	\end{equation}
	Clearly, \eqref{eq:comparison_2} allows a larger range of step-size compared to what \eqref{eq:comparison_1} does. This explains why ATC-DIGing performs better when graph is well-connected.
	For the analysis on ATC-DIGing, when $\delta$ is small, the left-hand-side of \eqref{eq:comparison_2} is at the order of $O(\delta^2\alpha)$ while the left-hand-side of \eqref{eq:comparison_3} is at the order of $O(\delta\alpha)$. Again, \eqref{eq:comparison_2} allows a larger range of step-size compared to what \eqref{eq:comparison_3} does. This is why we say a worse rate is derived in the current paper.
	
	Below we utilize \eqref{eq:comparison_1} and \eqref{eq:comparison_2} to show the complexity of DIGing and ATC-DIGing. In DIGing, we have
	\begin{equation}\label{eq:DIGing_bound}
	\lambda\leq1-\frac{(1-\delta)^2}{30n^{0.5}\bar{\kappa}^{1.5}},
	\end{equation}
	thus the iteration complexity to reach $\varepsilon$-accuracy is $O\left(\frac{n^{0.5}\bar{\kappa}^{1.5}}{(1-\delta)^2}\ln{\frac{1}{\varepsilon}}\right)$. In ATC-DIGing, the analogous of \eqref{eq:DIGing_bound} is
	\begin{equation}
	\lambda\leq1-\frac{(1-\delta)^2}{2(15n^{0.5}\bar{\kappa}^{1.5}\delta^2+1)}.
	\end{equation}
	Considering along with the limitation of gradient method [c.f. \eqref{eq:IGM_lambda}], we conclude that the iteration complexity of ATC-DIGing to reach $\varepsilon$-accuracy is
	$O\left(\max\left\{\frac{\delta^2n^{0.5}\bar{\kappa}^{1.5}+1}{(1-\delta)^2},\bar{\kappa}\right\}\ln{\frac{1}{\varepsilon}}\right)$. We omit the proof here since it is just a translation from convergence rate to computational complexity.
\end{document}